\documentclass[10pt,a4paper]{amsart}

\usepackage{amsfonts,amsmath,amssymb,graphics}

\input xy
\xyoption{all}

\newtheorem{thm}{Theorem}[section]
\newtheorem{lemma}[thm]{Lemma}
\newtheorem{cor}[thm]{Corollary}
\newtheorem{pro}[thm]{Proposition}
\newtheorem{example}[thm]{Example}
\newtheorem{remark}[thm]{Remark}
\newtheorem{ddef}[thm]{Definition}

\def\codim{\mathop{\rm codim}\nolimits}

\def\reg{\mathop{\rm reg}\nolimits}

\def\fin{\mathop{\rm end}}

\def\bideg{\mathop{\rm bideg}}
\def\cd{\mathop{\rm cd}\nolimits}

\def\codim{\mathop{\rm codim}}

\def\bideg{\mathop{\rm bideg}}

\def\C{{\mathcal C}}

\def\R{{\mathcal R}}

\def\im{{\mathfrak m}}

\def\ra{{\rightarrow}}

\def\fini{{$\quad\quad\Box$}}

\newcommand{\bd}{\begin{ddef}}
\newcommand{\ed}{\end{ddef}}
\newcommand{\bt}{\begin{thm}}
\newcommand{\et}{\end{thm}}
\newcommand{\bl}{\begin{lemma}}
\newcommand{\el}{\end{lemma}}
\newcommand{\bco}{\begin{cor}}
\newcommand{\eco}{\end{cor}}
\newcommand{\bp}{\begin{pro}}
\newcommand{\ep}{\end{pro}}
\newcommand{\bex}{\begin{example}}
\newcommand{\eex}{\end{example}}
\newcommand{\brm}{\begin{remark}}
\newcommand{\erm}{\end{remark}}
\newcommand{\bconj}{\begin{conj}}
\newcommand{\econj}{\end{conj}}

\newcommand{\fm}{\operatorname{\frak{M}}}

\newcommand{\fn}{\frak{n}}

\newcommand{\stab}{\operatorname{Stab}}

\newcommand{\beqn}{\begin{eqnarray*}}
\newcommand{\eeqn}{\end{eqnarray*}}
\newcommand{\beq}{\begin{eqnarray}}
\newcommand{\eeq}{\end{eqnarray}}
\newcommand{\been}{\begin{enumerate}}
\newcommand{\eeen}{\end{enumerate}}

\begin{document}

\author{Marc Chardin}
\address{Institut de Math\'ematiques de Jussieu\\
UPMC, Boite 247, 4, place Jussieu, F-75252 Paris Cedex, France\\}
\email{chardin@math.jussieu.fr}

\title[Regularity stabilization]{Regularity stabilization for the powers of graded $\fm$-primary ideals}

\maketitle

\begin{abstract}
This Note provides first a generalization of the stabilization result of Eisenbud and Ulrich 
for the regularity of
powers of a $\im$-primary ideal  to the case of ideals that are
not generated in a single degree (see Theorem \ref{stabregcoh}). We then partially extend our previous results expressing
this stabilization degree in term of the regularity of a specific graded strand of the Rees ring :
The natural extension of the statement holds at least if the stabilization index or the regularity is greater than the number
of variables. In any case a precise comparison is given (see Theorem \ref{SimpleStab}).

For simplicity, we do not
introduce a graded module as in \cite{EU}. It can be done along the same lines, but makes
the statements less transparent  (see  Remark \ref{modulecase} were we derive the key
point for such an extension).
\end{abstract}

\section{Stability and cohomology of a graded strand of the Rees ring}

Let $A:=k[x_1,\ldots ,x_n]$ be a standard graded polynomial ring over the field $k$, $\im :=(x_1,\ldots ,x_n)$ and $I=(f_1,\ldots ,f_m)$ be a  graded $\im$-primary ideal.

Set $d:=\min\{ \mu\ \vert\ \exists p, \ ( I_{\leq \mu})I^p =I^{p+1}\} $ and $I':= (I_{\leq d})$. By \cite{Ko}  there exists $b\geq 0$  such that
$$
\reg (I^t)=dt+b, \quad\forall t\gg0.
$$

Furthermore, 

\bp\label{t0}\cite[4.1]{Ch2}
Set $t_0:=\min\{ t\geq 1\ \vert\ \im^d I^t\subseteq I'\}$, then 

(i) $d=\min\{ \mu\ \vert \ (I_{\leq \mu})\ \hbox{is}\ \im\hbox{-primary} \}$.

(ii) The function $f(t):=\reg (I^t)-td$ is weakly decreasing for  $t\geq t_0$.

(iii)  One has,
$$
t_0\leq \max\left\{ 1,\left\lceil {{\reg (I')-d}\over{d+1}}\right\rceil\right\}
$$
and $\reg (I')\leq (d-1)n+1$.
\ep

In particular $t_0\leq \left\lceil {{(n-1)(d-1)}\over{d+1}}\right\rceil$ (unless $n=1$ or $d=1$).
 In Example 2.3 of [EU] $n=4$ and $d=5$ and this bound is sharp.\medskip

Similar related estimates concerning (ii) and (iii) are given in \cite{Be}.\medskip

We now want to derive results concerning the stabilization index : 
$$
\stab (I):=\min\{ t\ \vert\ \reg (I^s)=ds+b,\ \forall s\geq t\}.
$$
Let $(g_1,\ldots ,g_s)$ be generators of $I'=(I_{\leq d})$ with :\smallskip

(i) $\deg g_i =d$ for $1\leq i\leq m$ and  $\deg g_i\leq d$ for all $i$,\smallskip

(ii) setting  $J:=(g_1,\ldots ,g_m)$, $J+(I_{<d})$ is a reduction of $I$. \medskip

Condition (ii) is satisfied if $k$ is infinite, $m$ is at least $n-\ell ((I_{<d}))$,
where $\ell$ stands for the analytic spread, and the $g_i$'s are general elements 
in $I_d$.\medskip

Define, $B:=k[T_1,\ldots ,T_m]$, $\fn :=(T_1,\ldots ,T_m)$ and $S:=A[T_1,\ldots ,T_m] $.\par

We set $\bideg (T_i):=(0,1)$ and $\bideg (a):=(\deg (a),0)$ for $a\in A$. The natural inclusions
$\R_{J}\subset \R_{I'}\subset \R_I:=\oplus_{t\geq 0}I(d)^t=\oplus_{t\geq 0}I^t(td)$ and the onto map $S\ra \R_{J}$ makes $\R_I$ a bigraded $S$-module.
This module need not be finite over $S$. However $H^i_{\im}(\R_I )_{\mu ,*}$ is a 
finite graded $B$-module for every $\mu$ by \cite[2.1(ii)]{Ch1}, as $J+(I_{<d})$ is a reduction of $I$. \medskip

As $I$ is $\im$-primary, $H^1_\im (\R_I )=\oplus_{t\geq 1}(A/I^t)(td)$, $H^n_\im (\R_I )=\oplus_{t\geq 0}H^n_\im (A)(td)$, and
$H^i_{\im}(\R_I )=0$ for $i\not\in \{ 1,n\}$.\medskip

Denote by $N$ the finitely generated bigraded $S$-module $H^1_\im (\R_I )=\oplus_{t\geq 1}(A/I^t)(td)$.

\bl\label{stabilitytwodirections}
Let $t\geq 1$ and $\mu\geq 0$, then

(i) $N_{\mu ,t}=0\ \Rightarrow N_{\mu +1,t}=0$.

(ii) If $t\geq t_0$, then $N_{\mu ,t}=0\ \Rightarrow N_{\mu ,t+1}=0$.
\el

{\it Proof.} Claim (i) follows from the fact that $N_{*,t}=(A/I^t)(td)$ is cyclic
and generated in degree $-td\leq \mu$ and  (ii) is Proposition \ref{t0} (ii). \fini
\medskip

Part (ii) shows that $N_{\mu ,t}\not=0$ if $0\leq \mu <b$ and $t\geq t_0$.

Let 
$$
c:=\max_{t\geq 1}\{ \reg (I^t )-td\} =\reg_A (\R_I ).
$$
By Proposition \ref{t0}, $c$ is reached for a value of $t$ equal at most to $t_0$.

Claim (i) implies that $e_\mu :=\fin (N_{\mu ,*})$ is a weakly decreasing 
function of $\mu$ for $\mu\geq b$, which is equal to $-\infty$ for $\mu \geq c$. Hence,

\bp\label{stab1}
$$
\reg (I^t)=dt+b,\quad \forall t> e_b,
$$
and $\reg (I^t)>dt+b$ for any $t_0 \leq t\leq e_b$.
\ep

Notice that, as the examples in \cite[2.4]{EU} and \cite[Section 4]{Be} show, it may be that $N_{\mu ,t}=0$
for some $\mu \geq 0$ and $1\leq t<\min \{ e_\mu ,t_0\}$.  

\bl\label{1termSS}
There are natural bigraded isomorphisms
$$
H^1_{\im +\fn }(\R_I ) \simeq H^0_\fn (N) \simeq  H^0_{\im}(H^{1}_{\fn}(\R_I )) .
$$
\el

{\it Proof.} 
The two spectral sequences arising from the double complex $ \C^\bullet_\fn \C^\bullet_\im  \R_I$ abuting to $H^\bullet_{\fn +\im}(\R_I)$
have as first terms  $ \C^\bullet_\fn H^\bullet_\im  \R_I$ on one side. 
$H^0_\im  (\R_I )=0$, in particular, the component of cohomological degree $1$ in the abutment is
$H^0_\fn (H^1_{\im}(\R_I ))$. As $J$ has positive grade, $H^0_\fn (\R_I )=0$ and the other spectral sequence has 
only one term of cohomological degree $1$, which is $H^0_{\im}(H^{1}_{\fn}(\R_I ))$.
\fini
\medskip

This lemma shows that for $\mu \geq b$, $N_{\mu ,t}=H^0_\fn (N_{\mu ,t})\subseteq H^1_\fn ((\R_I )_{\mu ,*})_t$, which in turn shows that 
$$
N_{\mu ,t}=0,\quad \forall t\geq \reg_B ((\R_I )_{\mu ,*}).
$$
This implies the following,

\bt\label{stabregcoh} With notations as above,
$$
\reg (I^t)=dt+b,\quad \forall t\geq  \fin (H^1_\fn ((\R_I )_{b,*})+1,
$$
hence $\reg (I^t)=dt+b$, for $t\geq \reg_B ((\R_I )_{b,*})$. More precisely,
$$
\stab (I)=  \fin (H^0_{\im}(H^{1}_{\fn}(\R_I ))_{b,*})+1.
$$
\et

{\it Proof.} As $e_\mu :=\fin (N_{\mu ,*})$ is a weakly decreasing 
function of $\mu$ for $\mu\geq b$, the result follows from Proposition \ref{stab1} and the fact that
$$
e_b=\fin (H^0_{\im}(H^{1}_{\fn}(\R_I ))_{b,*})\leq \fin (H^1_\fn ((\R_I )_{b,*})\leq\reg_B ((\R_I )_{b,*})-1
$$
where the first equality is given by Lemma \ref{1termSS}. \fini

We will see in Proposition \ref{compcoh} that $ \fin (H^0_{\im}(H^{1}_{\fn}(\R_I ))_{b,*})= \fin (H^{1}_{\fn}(\R_I )_{b,*})$
if $J$ is $\im$-primary.

\brm\label{modulecase}
Lemma \ref{1termSS} extends to the case where a graded $A$-module $M$ 
such that $M/IM$ is $\im$-primary is involved, as
in the work of Eisenbud and Ulrich. In that case, the arguments in the proof of Lemma \ref{1termSS} show that
$$
H^1_{\im +\fn }(M\R_I ) \simeq  H^0_{\im}(H^{1}_{\fn}(M\R_I )) 
$$
 and $H^1_{\im +\fn }(M\R_I )_{*,t} \simeq  H^0_{\fn}(H^{1}_{\im}(M\R_I )) _{*,t}$ if $H^0_\im (M)=0$ or $t\geq 0$.
In that setting, $H^i_\im (MI^t)=H^i_\im (M)$ for $i\geq 2$ and one has an exact sequence
$$
0\to H^0_\im (I^tM)\to H^0_\im (M)\to M/I^tM\to H^1_\im (I^tM)\to H^1_\im (M)\to 0,
$$
which together with the variant of Lemma \ref{stabilitytwodirections} deduced from 
\cite{Ch2} provides an extension of \cite[1.1]{EU} to the unequal degree case.
\erm

\section{Finiteness of the regularity of graded strands of the Rees ring}

Notice first that $\reg_S (\R_I )<\infty$ if $\R_I$ is finite
over $S$, which in turn holds if $J$ is a reduction of $I$ ; hence in that case 
 $\reg_S ((\R_I )_{\mu ,*})\leq \reg_S (\R_I )<\infty $, for any $\mu$.

Let us also point out that $\reg_B (\R_I )=+\infty$ unless $I_{<d}=0$, as 
if not, then $(\R_I)_{*,t+1}$ cannot be generated over $B$ by $(\R_I)_{*,t}$, for any $t\geq 0$, for obvious degree
reason.

Next, the case of a complete intersection ideal will give us information. In that case one
can do an explicit computation, for instance as follows.

\brm
If $I$ is a complete intersection of $n$ forms of degree $d$, then one can compute $\dim (A/I^t)_{\mu +td}$ explicitly,  using 
\cite[2.3]{GVT}
$$
\dim_k (A/I^t)_{\mu +td}= \sum_{i=0}^{t-1}{{i+n-1}\choose{n-1}}\dim_k (A/I)_{\mu +(t-i)d},
$$
hence for $t\gg 0$ (more precisely if $\mu +(t+1)d>n(d-1)=\fin (A/I)$),
$$
\dim_k (A/I^t)_{\mu +td}= \sum_{j\geq 1}{{t-j+n-1}\choose{n-1}}\dim_k (A/I)_{\mu +jd}.
$$
It follows that $\dim_k (A/I^t)_{\mu +td}= D_\mu {{t+n-1}\choose{n-1}}+l.o.t.$ with $D_\mu :=\sum_{j\geq 1}\dim_k (A/I)_{\mu +jd}$. 

More generally, there exists a positive integer $C_\mu$ such that $\dim_k (A/I^t)_{\mu +td}= C_\mu {{t+m-1}\choose{m-1}}+l.o.t.$ if $I$ is a complete intersection 
of degrees $d_1,\ldots ,d_n$ with $d=d_1=\cdots =d_m>d_{m+1}\geq\cdots \geq d_n$, unless $\mu >\fin (A/I)-d$, in which case $\dim_k (A/I^t)_{\mu +td}=0$ for $t\geq 1$
\erm

The remark above illustrates the fact that the  finitely generated $S$-module $N_{\mu ,*}=H^1_\im (\R_I )_{\mu ,*}$ 
is of Krull dimension $d(\mu )\leq n-\ell (I_{<d})\leq n-\codim (I_{<d})$, and that this could be sharp. It in particular shows that $d(\mu )<n$ when $I_{<d}\not= 0$.
Thus the equality $d(\mu )=n$ implies that $I$ is generated in  degrees $\geq d$, which is equivalent to the fact that it admits a homogeneous reduction generated in degree $d$. 
 
 The following lemma provides a criteria for the finiteness of the regularity.
 
 \bl\label{reghp}
 Let $M$ be a graded $B$-module. Then the following are equivalent,
 
 (i) $\dim_k (M_\mu )<+\infty$ for some $\mu\geq \reg (M)$,
 
(ii) $\dim_k (M_\mu )<+\infty$ for all $\mu\geq \reg (M)$.

If these conditions hold, then there exist a polynomial $P$ of degree at most $m-1$ such that $\dim_k (M_\mu )=P(\mu )$ for
all $\mu >\reg (M)$.

\el

{\it Proof.} Let $H(\mu ):=\dim_k (M_\mu )$. (ii) implies (i) and to prove the converse, it suffices to show that $H(\mu +1)\geq H(\mu )$ for $\mu \geq \fin (H^0_\fn (M))$. 
Setting $\ell :=U_1T_1+\cdots +U_mT_m\in B\otimes_k k(U)$, $H(\mu )=\dim_{k(U)}M\otimes_k k(U)$ and 
$$
\ker (\!\!\xymatrix{M_\mu \otimes_k k(U)\ar^{\times \ell}[r]&M_{\mu +1} \otimes_k k(U)\\}\!\! )\subseteq H^0_\fn (M)_\mu ,
$$
by the Dedekind-Mertens Lemma (see e.g. \cite[1.7]{CJR}), which establishes our claim.

Notice that it suffices to show that  there exist a polynomial $P$ of degree at most $m-1$ such that $\dim_K (M_\nu )=P(\nu )$ for
all $\nu >\mu$ if $\mu\geq \reg (M)$, with $\mu\in {\bf Z}$. Let  $\mu\geq \reg (M)$ and $N:=M_{\geq \mu}$. One has $\reg (N)=\mu$, 
unless $N=0$, in which case $P=0$ satisfies our claim. It follows that $N$ is generated in degree $\mu$, hence is finitely generated. 
As $\reg (N)=\mu$, and $M$ and $N$ coincide in degree at least $\mu$, the conclusion follows.\fini

The useful consequence of Lemma \ref{reghp} for our study, is that any graded piece of the Rees algebra has
infinite regularity whenever the ideal $J$ has less than $n$ generators.

\bco
If $J$ has less than $n$ generators (i.e. $m<n$), then
$\reg_B (\R_I )_{\mu ,*}=+\infty$ for any $\mu\in {\bf Z}$. 
\eco

Notice that $m<n$  can only happen if $I_{<d}\not= 0$.\medskip

{\it Proof.}  By \cite[2.1(ii)]{Ch1}, the Hilbert function of
$N_{\mu ,*}$ is eventually a polynomial of degree at most $m-1$. As the dimension of $A_{\mu +td}$ 
over $k$ is a polynomial in $t$ of order $n-1$, it follows that it is also the case for  $(\R_I )_{\mu ,t}$
when $m<n$. From Lemma \ref{reghp}, we deduce the result as both (i) or (ii) are clear for $M:=(\R_I )_{\mu ,*}.$\fini

On the other hand, choosing $J$ generated by $n$ general elements in $I_d$ (if $k$ is infinite, else reducing to this case by
faithfully flat extension) the ideal $J+(I_{<d})$ is a reduction of $I$,  $N_{\mu ,t}$ is
finitely generated for any $\mu$ and $J$ is $\im$-primary.

\section{Cohomology and regularity of graded strands of the Rees ring}

In the context of the above two sections, the result below compares 
the cohomology of graded strands of $N=H^1_\im (\R_I )$ and of $\R_I$, in order to obtain
an estimate of the stabilization index in terms of regularity of a graded strands of 
the Rees ring.

\bp\label{compcoh}

If the ideal $J$ is $\im$-primary, then 

(i) $H^j_\fn (\R_I )_{* ,t}=H^0_\im (H^j_\fn (\R_I )_{* ,t})$ for any $j$ and $t\geq 0$.

(ii) For any $\mu >-n$ and any $j>0$, there exists a natural graded map of $B$-modules
$\psi_j^\mu : H^{j}_\fn (N_{\mu ,*})\to H^{j+1}_\fn (\R_I )_{\mu ,*}$, which is an 
isomorphism for $j>n$, onto for $j=n$, and an isomorphism in non negative degree
for any $j$.

(iii) For $\mu >-n$, let $d(\mu )$ be the Krull dimension of the support of the finitely generated graded
$B$-module $N_{\mu ,*}$. Then 
$$
\reg_B (N_{\mu ,*})\leq \max\{ \reg_B ((\R_I )_{\mu ,*} )-1,d(\mu ) \},
$$
$$
 \reg_B ((\R_I )_{\mu ,*})\leq \max\{ \reg_B (N_{\mu ,*})+1,\cd_{B_+}((\R_I )_{\mu ,*}) \}.
$$

(iv) $\cd_{B_+}((\R_I )_{\mu ,*})\leq n$ for $\mu >-n$. 

\ep

This in particular shows that $ \reg_B (N_{\mu ,*})= \reg_B ((\R_I )_{\mu ,*} )-1$ if either 

(1) $ \reg_B (N_{\mu ,*})\geq n$ and $I_{<d}\not= 0$; or 

(2) $\max\{  \reg_B (N_{\mu ,*}), \reg_B ((\R_I )_{\mu ,*})\} >n$.

{\it Proof.} First notice that after inverting any of the variables $X_i$, one has
graded (but in general not bigraded) identifications
$$
(\R_{J})_{X_i} =  (\R_J)_{X_i} =  (\R_I)_{X_i} = \oplus_{t\geq 0}A_{X_i}T^t =A_{X_i}[T],
$$
as $JA_{X_i}=JA_{X_i}=IA_{X_i}=A_{X_i}$. In the above identification, the positive part $\fn$
of $S$ is mapped to the ideal generated by $TJA_{X_i}$ which is the ideal $(T)=\oplus_{t> 0}A_{X_i}T^t$.

Hence, the spectral sequence with first terms $E^{i,j}_1= \C^i_\im  H^j_\fn \R_I=H^j_\fn \C^i_\im  \R_I$
satisfies $E^{i,j}_1=0$ unless $i=0$, or $j=1$. Also notice that $(E^{i,1})_{*,t}=0$ for
$t\geq 0$ and $i>0$ by virtue of the identity $H^1_{(T)} ( A_{X_i}[T])=T^{-1}A_{X_i}[T^{-1}]$, proving claim (i).

Furthermore, its second terms $E^{i,j}_2=H^i_\im (H^j_\fn   \R_I)$ vanish 
for $i>n$, hence 
\[
'E^{i,j}_2=\left\{
\begin{array}{cl}
H^j_\fn (\R_I)&\hbox{if $i=0$ and $j\not= 1$}\\
H^i_\im (H^1_\fn (\R_I ))& \hbox{if $0\leq i\leq n$ and $j= 1$}\\
0& \hbox{else}\\
\end{array}
\right.
\]

It shows that $(E^{i,j}_2)_{\mu ,t}=(E^{i,j}_\infty)_{\mu ,t}=0$ for $i\not= 0$, if 
either $t\geq 0$ or $j>n$, which implies that  $(H^j_\fn (\R_I))_{\mu ,t}\simeq (E^{0,j}_2)_{\mu ,t}\simeq (E^{0,j}_\infty)_{\mu ,t}$ if 
either $t\geq 0$ or $j>n$.

On the other hand, $H^1_\im (\R_I )=\oplus_{t\geq 1} (R/I^t) (td)$,  $H^n_\im (\R_I )=\oplus_{t\geq 0} H^n_\im (R) (td)$, and $H^i_\im (\R_I )=0$
for $i\not\in \{ 1,n\}$.

The second spectral sequence has second terms graded components :
\[
(''E^{i,j}_2)_{\mu ,*}=\left\{
\begin{array}{cl}
H^j_\fn (H^1_\im (\R_I )_{\mu ,*})&\hbox{if $i=1$}\\
H^j_\fn (H^n_\im (\R_I )_{\mu ,*})& \hbox{if $i=n$}\\
0& \hbox{else}\\
\end{array}
\right.
\]
showing that the $(''E^{i,j}_2)_{\mu ,*}=0$  for $i\not= 1$ and $\mu >-n$.
This in turn shows that, if $\mu >-n$, there is a long exact sequence 
$$
\cdots\to H^{j-1}_\im (H^1_\fn (\R_I )_{\mu ,*})\to H^{j-1}_\fn (N_{\mu ,*})\to H^j_\fn (\R_I )\to H^j_\im (H^1_\fn (\R_I )_{\mu ,*})\to \cdots,
$$
where the map $H^j_\fn (\R_I )\to H^j_\im (H^1_\fn (\R_I )_{\mu ,*})$ is $'d^{0,j}_j$ for $2\leq j\leq n$ and 0 else.

$H^{n+1}_\fn (\R_I )_{\mu ,*}=0$, and 
$$
H^j_\fn (\R_I )_{\mu ,t}\simeq H^{j-1}_\fn (H^1_\im (\R_I )_{\mu ,*})_t,\quad \forall j, t\geq 0, \mu >-n.
$$
for $j\geq 2$ and $H^1_\im (H^1_\fn (\R_I ))_{\mu ,*}=0$. In particular, $H^j_\fn (\R_I )_{\mu ,t-j}=0$ if $j\geq 2$, $t\geq n$ and $\mu >-n$.

Now, $H^j_\fn (\R_I )_{\mu ,t-j}=0$, $\forall j$ implies $H^j_\fn (\R_I )_{\mu ,t+1-j}=0$, $\forall j$, as $H^j_\fn (\R_I )_{\mu ,s}=(H^j_\fn (\R_I )_{\mu ,*})_{s}$ for any
$\mu$ and $s$.

This shows that for $t\geq n$,  $H^1_\im (\R_I )_{\mu ,t-1}=0$ implies $H^1_\im (\R_I )_{\mu ,t}=0$. \fini

\section{A simple stability result}

Out of the results of section 1 and 3, we can deduce the following result that provides a
quite sharp estimate of the stabilization index in terms of the regularity of graded strands 
of the Rees ring for non equigenerated ideals.

\bt\label{SimpleStab}
Let $A$ be a standard graded polynomial ring in $n$ variables over an infinite field, $\im:=A_+$ and $I$ be a $\im$-primary graded  ideal. 

Set $d:=\min\{ \mu\ \vert \ (I_{\leq \mu})\ \hbox{is}\ \im\hbox{-primary} \}$ and $t_0:=\min\{ t\geq 1\ \vert\ \im^d I^t\subseteq (I_{\leq d})\}$. 

Let $J=(g_1,\ldots ,g_n)\subseteq I$ be generated by $n$ general elements in $I_d$. 

Consider the standard bigraded ring $S:=A[T_1,\ldots ,T_n]$, $B:=k[T_1,\ldots ,T_n]$, $\fn :=(T_1,\ldots ,T_n)$ and
denote by $\R_J= S/(\{ g_iT_j-g_jT_i\ \vert\ 1\leq i<j\leq n\} )$ the Rees ring of the complete intersection ideal $J(d)$. 

The inclusion $\R_J\subset \R_I:=\oplus_t I(d)^t$ makes $\R_I$ a bigraded module over  $S$.

Let $b$ be defined by $\reg (I^t)=dt+b$ for $t\gg 0$. Then
$$
\stab (I)= \fin (H^1_{\fn}(\R_I)_{b,*})+1
$$
and
$$
\stab (I)\leq  \reg_B ((\R_I )_{b ,*})\leq \max\{ \stab (I),n \}.
$$

In particular, $\stab (I)=\reg_B ((\R_I )_{b ,*})$ if $\stab (I)\geq n$ or $\reg_B ((\R_I )_{b ,*})>n$.
\et

Notice that if $I$ has no generator of degree less than $d$, then  $J$ is a minimal reduction of 
$I$. Also recall that $t_0\leq n-1$  by Proposition \ref{t0}.

\end{document}